\documentclass{amsart}

\usepackage[utf8]{inputenc}
\usepackage{amsmath,amssymb,amsthm}
\usepackage{latexsym}
\usepackage{hyperref,mathtools}
\hypersetup{colorlinks=true}
\usepackage{tikz}
\usepackage{graphics}
\usepackage{dsfont}
\usepackage{bbm}
\usepackage{siunitx}

\newtheorem*{thm*}{Theorem}

\newcommand{\Z}{\mathbb{Z}}
\newcommand{\Prob}{\mathbb{P}}

\newcommand{\vL}{\vec{\mathbb{L}}_2}

\newcommand{\vLt}{\vec{\mathbb{L}}_3}

\newcommand{\Lva}{\vec{\mathbb{L}}_3}

\newcommand{\Laltd}{\vec{\mathbb{L}}_2^\mathrm{alt}}
\title[Update on lower bounds for oriented percolation]{An update on lower bounds for the critical values of oriented percolation models}
\author{Olivier Couronné}
\address{Universit\'e Paris Nanterre, Modal'X, FP2M, CNRS FR 2036, 200 avenue de la R\'epublique 92000 Nanterre, France.}
\email{olivier.couronne@parisnanterre.fr}
\thanks{The author is supported by the Labex MME-DII funded by ANR,
	reference ANR-11-LBX-0023-01, and this research has been conducted within the FP2M federation (CNRS FR 2036)}

\begin{document}
	
	\maketitle
	\begin{abstract}
		We obtain new lower bounds on the critical points for various models of oriented percolation. The method is to provide a stochastic domination of the percolation processes by multitype Galton-Watson trees. This can be apply to the classical bond and site oriented percolation on $\Z^2$, but also on other lattices such as inhomogeneous ones, and on dimension three.
	\end{abstract}
	\section{Introduction}
	We develop a new method for studying oriented percolation processes, method which provides lower bounds for the threshold of these processes. A classical way to obtain a lower bound in percolation is via a Peierl argument, where one enumerates all possible (self-avoiding when we can) paths from the origin $O$ of length $n$, and the threshold is larger than the geometric progression in $n$ of this sequence. In this article, we do not consider all open paths (using only
	rightward and upward edges), but all paths such that there is no other path below that attains the same vertex. Such paths will be called valid paths. The process obtained has less open edges but the cluster of the origin is the same than that of the classical model.

	In 2006, \cite{belitsky} provided a lower bound for the threshold of the classical bond oriented percolation on $\Z^2$. Their method was based on the asymptotic and renormalized position of the rightmost infected particle when starting with all particles on the left of the origin infected. 
	There are pros and cons for their method and ours, but our method can readily be applied in other  two-dimensional oriented situations, detailed in section \ref{sec:otherTwoDim}, and also in three dimensions, as explained in section~\ref{sec:threeDim}.
	
	The main results of this paper are the following:
	\vskip -3mm
	\begin{center}
		\begin{tabular}{ |l|l|l|l| } 
			\hline
			Model & New lower bound & Previous one  &Estimate\\ 
			\hline
			site percolation on $\vL$ & $0.6967$ & $0.6882$ \cite{Gray} & $0.7055$ \cite{Wang_2013}\\ 
			\hline
			bond percolation on $\vL$ &$0.636893$ & $0.63328$  \cite{belitsky} & $0.6447$ \cite{Wang_2013}\\ 
			\hline
			site percolation on $\Laltd$& $0.525$&$0.5$ \cite{Pearce} & $0.535$ \cite{Pearce}
			\\
			\hline
			bond percolation on $\Laltd$&$0.4022$ &$1/3$ (Peierl) & -\\
			\hline
			site percolation on $\vLt$ & $0.41507$ & $1/3$ (Peierl) & $0.43531$ \cite{Wang_2013}\\
			\hline
			bond percolation on $\vLt$ & $0.360684$ & $1/3$ (Peierl) & $0.38222$ \cite{Wang_2013}\\
			\hline 
		\end{tabular}
	\end{center}
	We postpone the results for inhomogeneous models in dimension two to section~\ref{sec:otherTwoDim}.
	
	The paper is organized as follows. Section~\ref{oriented2d} is devoted to the classical oriented percolation on $\Z^2$. We describe, in the context of the standard oriented percolation on $\Z^2$, the multitype Galton-Watson process in section~\ref{sec:automaton}. 
	We then prove in section~\ref{sec:domination} that if the Galton-Watson process is subcritical, then so is the percolation process.
	We give in section~\ref{sec:comparaison} some elements that may help the reader to compare, in terms of efficiency, our method to the one used in~\cite{belitsky}.
	In section~\ref{sec:otherTwoDim}, we give lower bounds for the threshold on some inhomogeneous oriented percolation models. In section~\ref{sec:threeDim}, we show how to adapt the states in the three dimensional setting.
	
	\section{Oriented percolation on $\Z^2$}
	\label{oriented2d}
	\subsection{The model}
	The set of vertices is $\Z^2$. From each $(x, y)\in \Z^2$, there is an oriented bond to $(x+1, y)$ and another one to $(x, y+1)$. Below is the classical representation of this graph, together with an equivalent one we shall use in the sequel:
	
	
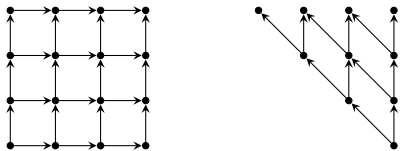
\begin{center}
		\begin{tikzpicture}[scale=0.6, ->]
			\draw[>=stealth, shorten >= 1.5pt](0,0)--(1,0);
			\draw[>=stealth, shorten >= 1.5pt](0,0)--(0,1);
			\draw[>=stealth, shorten >= 1.5pt](1,0)--(2,0);
			\draw[>=stealth, shorten >= 1.5pt](1,0)--(1,1);
			\draw[>=stealth, shorten >= 1.5pt](2,0)--(3,0);
			\draw[>=stealth, shorten >= 1.5pt](2,0)--(2,1);
			\draw[>=stealth, shorten >= 1.5pt](3,0)--(3,1);
			\filldraw (3,0) circle (2pt);
			\filldraw (2,0) circle (2pt);
			\filldraw (1,0) circle (2pt);
			\filldraw (0,0) circle (2pt);
			\begin{scope}[shift={(0,1)}]
				\draw[>=stealth, shorten >= 1.5pt](0,0)--(1,0);
				\draw[>=stealth, shorten >= 1.5pt](0,0)--(0,1);
				\draw[>=stealth, shorten >= 1.5pt](1,0)--(2,0);
				\draw[>=stealth, shorten >= 1.5pt](1,0)--(1,1);
				\draw[>=stealth, shorten >= 1.5pt](2,0)--(3,0);
				\draw[>=stealth, shorten >= 1.5pt](2,0)--(2,1);
				\draw[>=stealth, shorten >= 1.5pt](3,0)--(3,1);
				\filldraw (3,0) circle (2pt);
				\filldraw (2,0) circle (2pt);
				\filldraw (1,0) circle (2pt);
				\filldraw (0,0) circle (2pt);
			\end{scope}
			\begin{scope}[shift={(0,2)}]
				\draw[>=stealth, shorten >= 1.5pt](0,0)--(1,0);
				\draw[>=stealth, shorten >= 1.5pt](0,0)--(0,1);
				\draw[>=stealth, shorten >= 1.5pt](1,0)--(2,0);
				\draw[>=stealth, shorten >= 1.5pt](1,0)--(1,1);
				\draw[>=stealth, shorten >= 1.5pt](2,0)--(3,0);
				\draw[>=stealth, shorten >= 1.5pt](2,0)--(2,1);
				\draw[>=stealth, shorten >= 1.5pt](3,0)--(3,1);
				\filldraw (3,0) circle (2pt);
				\filldraw (2,0) circle (2pt);
				\filldraw (1,0) circle (2pt);
				\filldraw (0,0) circle (2pt);
			\end{scope}
			\begin{scope}[shift={(0,3)}]
				\filldraw (3,0) circle (2pt);
				\filldraw (2,0) circle (2pt);
				\filldraw (1,0) circle (2pt);
				\filldraw (0,0) circle (2pt);
				\draw[>=stealth, shorten >= 1.5pt](0,0)--(1,0);
				\draw[>=stealth, shorten >= 1.5pt](1,0)--(2,0);
				\draw[>=stealth, shorten >= 1.5pt](2,0)--(3,0);
			\end{scope}
			
			\begin{scope}[shift={(8.5,0)}]
				\filldraw (0,0) circle (2pt);
				\filldraw (-1,1) circle (2pt);
				\filldraw (0,1) circle (2pt);
				\filldraw (-1,2) circle (2pt);
				\filldraw (0,2) circle (2pt);
				\filldraw (-2,2) circle (2pt);
				\filldraw (-1,3) circle (2pt);
				\filldraw (0,3) circle (2pt);
				\filldraw (-2,3) circle (2pt);
				\filldraw (-3,3) circle (2pt);
				\draw[>=stealth, shorten >= 1.5pt](0,0)--(0,1);
				\draw[>=stealth, shorten >= 1.5pt](0,0)--(-1,1);
				\begin{scope}[shift={(0,1)}]
					\draw[>=stealth, shorten >= 1.5pt](0,0)--(0,1);
					\draw[>=stealth, shorten >= 1.5pt](0,0)--(-1,1);
				\end{scope}
				\begin{scope}[shift={(-1,1)}]
					\draw[>=stealth, shorten >= 1.5pt](0,0)--(0,1);
					\draw[>=stealth, shorten >= 1.5pt](0,0)--(-1,1);
				\end{scope}
				\begin{scope}[shift={(0,2)}]
					\draw[>=stealth, shorten >= 1.5pt](0,0)--(0,1);
					\draw[>=stealth, shorten >= 1.5pt](0,0)--(-1,1);
				\end{scope}
				\begin{scope}[shift={(-1,2)}]
					\draw[>=stealth, shorten >= 1.5pt](0,0)--(0,1);
					\draw[>=stealth, shorten >= 1.5pt](0,0)--(-1,1);
				\end{scope}
				\begin{scope}[shift={(-2,2)}]
					\draw[>=stealth, shorten >= 1.5pt](0,0)--(0,1);
					\draw[>=stealth, shorten >= 1.5pt](0,0)--(-1,1);
				\end{scope}
			\end{scope}
		\end{tikzpicture}
	\end{center}
	The graph is denoted $\vL$. The {\it height} of a vertex $(x, y)$ is $x+y$.
	The {\it successors} of a vertex are the endpoints of the oriented edges starting from the vertex. Here a vertex $(x, y)$ has $(x+1, y)$ and $(x, y+1)$ for successors.
	
	Fix a parameter $p\in [0, 1]$. 
	We can consider {\it site} or {\it bond} percolation on the graph:
	\begin{itemize}
		\item For site percolation, each vertex is open with probability $p$, closed with probability $1-p$, independently of each other.
		\item For bond percolation, each bond
		is open with probability $p$, closed with probability $1-p$, independently of each other.
	\end{itemize}
	We write $\Prob$ for the studied law on the graph. The critical value for $p$ is then
	\[p_c=\inf\{p>0, \Prob(\mbox{there exists an infinite open path from }0)>0\}.\]
	See \cite{Grimmett} for a classical reference on percolation.
	
	We say that a vertex is {\it occupied} if there is an open oriented path from the origin to the vertex, and {\it vacant} otherwise. As there is at most two vertices directly connected to each vertex, by a Peierl argument, we see that $p_c\geq 0.5$. This lower bound is not satisfactory, and this is due to the fact that a lot of counted paths have the same endpoints, so there is redundancy. 
	
	In order to keep the description simple, we now restrict ourselves to the case of oriented bond percolation.
	We propose another process which has the same occupied vertices, but for edges a subset of the open edges. We say that an edge is {\it good} if it is open, starts from an occupied vertex denoted $(x, y)$, and either
	\begin{itemize}
		\item is vertical.
		\item is horizontal and there is no open path from the origin to $(x+1, y)$ passing through $(x+1, y-1)$.
	\end{itemize}
	
	
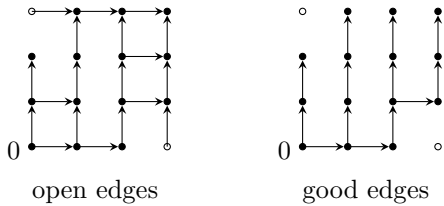
\begin{center}
		\begin{tikzpicture}[scale=0.6, ->]
			\node at(-0.4, -0.08) {$0$};
			\node at(1.4, -1) {open edges};
			\draw[>=stealth, shorten >= 1.5pt](0,0)--(1,0);
			\draw[>=stealth, shorten >= 1.5pt](0,0)--(0,1);
			\draw[>=stealth, shorten >= 1.5pt](1,0)--(2,0);
			\draw[>=stealth, shorten >= 1.5pt](1,0)--(1,1);
			\draw[>=stealth, shorten >= 1.5pt](2,0)--(2,1);
			\draw[>=stealth, shorten >=1.5pt] (3,0)--(3,1);
			\draw (3,0) circle (2pt);
			\filldraw (2,0) circle (2pt);
			\filldraw (1,0) circle (2pt);
			\filldraw (0,0) circle (2pt);
			\begin{scope}[shift={(0,1)}]
				\draw[>=stealth, shorten >= 1.5pt](0,0)--(1,0);
				\draw[>=stealth, shorten >= 1.5pt](0,0)--(0,1);
				\draw[>=stealth, shorten >= 1.5pt](1,0)--(1,1);
				\draw[>=stealth, shorten >= 1.5pt](2,0)--(3,0);
				\draw[>=stealth, shorten >= 1.5pt](2,0)--(2,1);
				\draw[>=stealth, shorten >= 1.5pt](3,0)--(3,1);
				\filldraw (3,0) circle (2pt);
				\filldraw (2,0) circle (2pt);
				\filldraw (1,0) circle (2pt);
				\filldraw (0,0) circle (2pt);
			\end{scope}
			\begin{scope}[shift={(0,2)}]
				\draw[>=stealth, shorten >= 1.5pt](1,0)--(1,1);
				\draw[>=stealth, shorten >= 1.5pt](2,0)--(3,0);
				\draw[>=stealth, shorten >= 1.5pt](2,0)--(2,1);
				\draw[>=stealth, shorten >= 1.5pt](3,0)--(3,1);
				\filldraw (3,0) circle (2pt);
				\filldraw (2,0) circle (2pt);
				\filldraw (1,0) circle (2pt);
				\filldraw (0,0) circle (2pt);
			\end{scope}
			\begin{scope}[shift={(0,3)}]
				\filldraw (3,0) circle (2pt);
				\filldraw (2,0) circle (2pt);
				\filldraw (1,0) circle (2pt);
				\draw (0,0) circle (2pt);
				\draw[>=stealth, shorten >= 1.5pt](0,0)--(1,0);
				\draw[>=stealth, shorten >= 1.5pt](1,0)--(2,0);
				\draw[>=stealth, shorten >= 1.5pt](2,0)--(3,0);
			\end{scope}
			
			\begin{scope}[shift={(6,0)}]
				\node at(1.4, -1) {good edges};
				\node at(-0.4, -0.08) {$0$};
				\draw[>=stealth, shorten >= 1.5pt](0,0)--(1,0);
				\draw[>=stealth, shorten >= 1.5pt](0,0)--(0,1);
				\draw[>=stealth, shorten >= 1.5pt](1,0)--(2,0);
				\draw[>=stealth, shorten >= 1.5pt](1,0)--(1,1);
				\draw[>=stealth, shorten >= 1.5pt](2,0)--(2,1);
				\draw (3,0) circle (2pt);
				\filldraw (2,0) circle (2pt);
				\filldraw (1,0) circle (2pt);
				\filldraw (0,0) circle (2pt);
				\begin{scope}[shift={(0,1)}]
					\draw[>=stealth, shorten >= 1.5pt](0,0)--(0,1);
					\draw[>=stealth, shorten >= 1.5pt](1,0)--(1,1);
					\draw[>=stealth, shorten >= 1.5pt](2,0)--(3,0);
					\draw[>=stealth, shorten >= 1.5pt](2,0)--(2,1);
					\draw[>=stealth, shorten >= 1.5pt](3,0)--(3,1);
					\filldraw (3,0) circle (2pt);
					\filldraw (2,0) circle (2pt);
					\filldraw (1,0) circle (2pt);
					\filldraw (0,0) circle (2pt);
				\end{scope}
				\begin{scope}[shift={(0,2)}]
					\draw[>=stealth, shorten >= 1.5pt](1,0)--(1,1);
					\draw[>=stealth, shorten >= 1.5pt](2,0)--(2,1);
					\draw[>=stealth, shorten >= 1.5pt](3,0)--(3,1);
					\filldraw (3,0) circle (2pt);
					\filldraw (2,0) circle (2pt);
					\filldraw (1,0) circle (2pt);
					\filldraw (0,0) circle (2pt);
				\end{scope}
				\begin{scope}[shift={(0,3)}]
					\filldraw (3,0) circle (2pt);
					\filldraw (2,0) circle (2pt);
					\filldraw (1,0) circle (2pt);
					\draw (0,0) circle (2pt);
				\end{scope}
			\end{scope}
		\end{tikzpicture}
	\end{center}
	
	With these rules, there is only one path of good edges connecting the origin to an occupied vertex. This will be used to provide a new lower bound on the critical point.
	\subsection{The automaton for the oriented  percolation on $\Z^2$.}
	\label{sec:automaton}
	Fix an integer $k\geq2$. A state is a sequence of length $k$ of $0$ and $1$, starting with a $1$. Essentially, it corresponds to the occupancy (associated to $1$) and vacancy (associated to $0$) of an interval of $k$ vertices of the same height, with the leftmost vertex being occupied.
	
	To define the children of a state, we consider it as an horizontal succession of occupied and vacant vertices on the second representation of the graph, together with the edges starting from these vertices, and their successors (we have here  $k+1$ successors).
	See the following picture for an example with $k=4$:
	
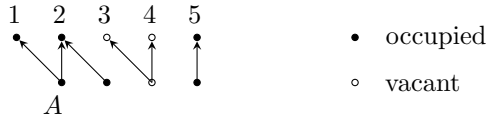
\begin{center}
		\begin{tikzpicture}[scale=0.6, ->]
			\draw[>=stealth, shorten >= 1.5pt](0,0)--(-1,1);
			\draw[>=stealth, shorten >= 1.5pt](0,0)--(0,1);
			\draw[>=stealth, shorten >= 1.5pt](1,0)--(0,1);
			\draw[>=stealth, shorten >= 1.5pt](2,0)--(2,1);
			\draw[>=stealth, shorten >= 1.5pt](2,0)--(1,1);
			\draw[>=stealth, shorten >= 1.5pt](3,0)--(3,1);
			\filldraw (3,0) circle (2pt);
			\draw (2,0) circle (2pt);
			\filldraw (1,0) circle (2pt);
			\filldraw (0,0) circle (2pt);
			\filldraw (3,1) circle (2pt);
			\draw (2,1) circle (2pt);
			\draw (1,1) circle (2pt);
			\filldraw (0,1) circle (2pt);
			\filldraw (-1,1) circle (2pt);
			\node at (-0.2,-0.5) {$A$};
			\node at (-1.05,1.5) {$1$};
			\node at (-0.05,1.5) {$2$};
			\node at (0.95,1.5) {$3$};
			\node at (1.95,1.5) {$4$};
			\node at (2.95,1.5) {$5$};
			\draw (6.5,0) circle (2pt);
			\filldraw (6.5,1) circle (2pt);
			\node at (8.3,1) {occupied};
			\node at (8,0) {vacant};
		\end{tikzpicture}
	\end{center}
	The children are the intervals of length $k$, starting with a vertex infected by the leftmost vertex of the state, {\it via} a good edge. In the figure above, the vertex $1$ is indeed infected via a good edge, but concerning the vertex $2$, the infection from $A$ does not use a good edge because of the other open edge infecting vertex $2$. In this case, the state has one child, that is $(1100)$. Below we list several examples:
	
	
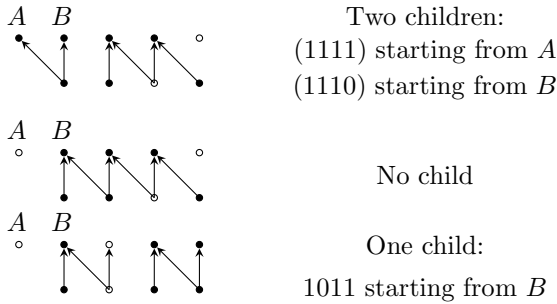
\begin{tikzpicture}[scale=0.6, ->]
		\draw[>=stealth, shorten >= 1.5pt](0,0)--(-1,1);
		\draw[>=stealth, shorten >= 1.5pt](0,0)--(0,1);
		\draw[>=stealth, shorten >= 1.5pt](1,0)--(1,1);
		\draw[>=stealth, shorten >= 1.5pt](2,0)--(2,1);
		\draw[>=stealth, shorten >= 1.5pt](2,0)--(1,1);
		\draw[>=stealth, shorten >= 1.5pt](3,0)--(2,1);
		\filldraw (3,0) circle (2pt);
		\draw (2,0) circle (2pt);
		\filldraw (1,0) circle (2pt);
		\filldraw (0,0) circle (2pt);
		\draw (3,1) circle (2pt);
		\filldraw (2,1) circle (2pt);
		\filldraw (1,1) circle (2pt);
		\filldraw (0,1) circle (2pt);
		\filldraw (-1,1) circle (2pt);
		\node at (8,1.5) {Two children:};
		\node at (8,0.7) {$(1111)$ starting from $A$};
		\node at (8,-0.1) {$(1110)$ starting from $B$};
		\node at (-1.05,1.5) {$A$};
		\node at (-0.05,1.5) {$B$};
	\end{tikzpicture}
	
	\begin{tikzpicture}[scale=0.6, ->]
		\draw[>=stealth, shorten >= 1.5pt](1,0)--(0,1);
		\draw[>=stealth, shorten >= 1.5pt](0,0)--(0,1);
		\draw[>=stealth, shorten >= 1.5pt](1,0)--(1,1);
		\draw[>=stealth, shorten >= 1.5pt](2,0)--(2,1);
		\draw[>=stealth, shorten >= 1.5pt](2,0)--(1,1);
		\draw[>=stealth, shorten >= 1.5pt](3,0)--(2,1);
		\filldraw (3,0) circle (2pt);
		\draw (2,0) circle (2pt);
		\filldraw (1,0) circle (2pt);
		\filldraw (0,0) circle (2pt);
		\draw (3,1) circle (2pt);
		\filldraw (2,1) circle (2pt);
		\filldraw (1,1) circle (2pt);
		\filldraw (0,1) circle (2pt);
		\draw (-1,1) circle (2pt);
		\node at (8,0.5) {No child};
		\node at (-1.05,1.5) {$A$};
		\node at (-0.05,1.5) {$B$};
	\end{tikzpicture}

	\begin{tikzpicture}[scale=0.6, ->]
		\draw[>=stealth, shorten >= 1.5pt](1,0)--(0,1);
		\draw[>=stealth, shorten >= 1.5pt](0,0)--(0,1);
		\draw[>=stealth, shorten >= 1.5pt](1,0)--(1,1);
		\draw[>=stealth, shorten >= 1.5pt](2,0)--(2,1);
		
		\draw[>=stealth, shorten >= 1.5pt](3,0)--(2,1);
		\draw[>=stealth, shorten >= 1.5pt](3,0)--(3,1);
		\filldraw (3,0) circle (2pt);
		\filldraw (2,0) circle (2pt);
		\draw (1,0) circle (2pt);
		\filldraw (0,0) circle (2pt);
		\filldraw (3,1) circle (2pt);
		\filldraw (2,1) circle (2pt);
		\draw (1,1) circle (2pt);
		\filldraw (0,1) circle (2pt);
		\draw (-1,1) circle (2pt);
		\node at (8,1) {One child:};
		\node at (8,0) {$1011$ starting from $B$};
		\node at (-1.05,1.5) {$A$};
		\node at (-0.05,1.5) {$B$};
	\end{tikzpicture}
	
	\noindent
	If the state has a child starting from $A$, we say that it is its {\it up-child}. If the state has a child starting from $B$, we say that it is its {\it right-child}.
	
	Hence the law of the offspring of a state is given by the state and the law of the open edges starting from the vertices of the state. We define in that way a multitype Galton-Watson process, with its root state composed of one $1$ followed by $0$'s ($(1000)$ for $k=4$). We refer to \cite{athreya} for a study of these processes.
	
	To know if the Galton-Watson process is either subcritical or supercritical, we enumerate all the states and consider the transition matrix, where the term at $(i, j)$ is the expectation of the number of states $j$ obtained from state $i$. If the largest eigenvalue of this matrix is strictly larger than $1$, then the process is supercritical. Otherwise it is subcritical.
	
	We recall that each state has at most two children, one up-child and one right-child. We consider the Galton-Watson process as a subset of a regular and labelled $2$-tree, going upwards, with the state up-child associated to the tree child on the left, and the state right-child associated to the tree child on the right. 
	We call $T$ the regular $2$-tree, without the labels.
	\subsection{Stochastic domination}\label{sec:domination}
	In this section we explicit the stochastic domination of the percolation process by the multitype Galton-Watson tree described in the preceding section. For an integer $n$, let $L_n$ be the set of vertices of the form $(x, n-x)$. In the following sum, the paths $\gamma$ are paths on $\Z^2$ using only up and right steps. We say that a path is {\it valid} if it uses only good vertices.
	
	\begin{equation}
		\Prob(O\rightarrow L_n)\leq \sum_{\gamma: O\rightarrow L_n}  \Prob(\gamma \text{ is valid}) 
	\end{equation}
	
	For a given path $\gamma$, we associate a coupling between the percolation and the Galton-Watson processes. We will always start for the Galton-Watson process with the state $s_0=(10\ldots0)$.
	To the path $\gamma$ corresponds a path, also denoted $\gamma$, on the regular $2$-tree $T$. 
	
	We build a sequence of states $(s_m)_{0\leq m\leq n}$ as follows: while the state $s_m$ is defined, the edges used to get the children of the state $s_m$ are in the same configuration as the edges starting from the diagonal segment of length $k$ and beginning at $\gamma_m$ (the $(m+1)$th vertex of $\gamma$). 
	This defines the (possibly null) up-child and right-child of $s_m$. Now
	\begin{itemize}
		\item If $\gamma_{m+1}$ is the upward successor of $\gamma_m$, and $\gamma_{m+1}$ is occupied, then $s_{m+1}$ is the up-child of $s_m$.
		\item If $\gamma_{m+1}$ is the rightward successor of $\gamma_m$, and $s_m$ has a right-child, then $s_{m+1}$ is the right-child of $s_m$.
	\end{itemize}
	In particular, if $\gamma_{m+1}$ is vacant, then $s_m$ has no child and the coupling stops here for this branch.
	
	For all the offsprings not yet considered, the Galton-Watson process evolves according to the transitions previously defined.
	
	The path $\gamma$ is still fixed. One can see, by recurrence, that for each $m\geq 0$ such that $s_m$ exists, the state $s_m$ is less than the status of the corresponding vertices, that is if there is a $1$ in the state, then the corresponding vertex is occupied. From this we deduce that if the edge from  
	$\gamma_m$ to $\gamma_{m+1}$ is good, then the state $s_m$ has a child in the corresponding direction. 
	Recall that an vertical edge from an occupied vertex is good, so this remark concerns only rightwards edges. 
	
	Finally, with this coupling, we have obtained that if $\gamma$ is valid,
	then $\gamma_n$  is alive in the Galton-Watson process. Below we write $T_{n, i}$ for the $i$th vertex of height $n$ on the tree $T$ (there is $2^n$ of them for a given $n$). Summing over all the paths $\gamma$ from $O$ to $L_n$, this gives
	
	\[\Prob(O\rightarrow L_n)\leq \sum_{i=1}^{2^n}\Prob( T_{n,i}\text{ is alive}).\]
	But when the Galton-Watson process is subcritical, the term on the right, being the expectation of the number of particles alive at time $n$, converges to $0$ (see \cite{athreya}). We have proved that when the multitype Galton-Watson process is subcritical, so is the oriented percolation process.
	
	\subsection{Some remarks concerning the implementation}
	For a given $k$, the space of the states is of cardinal $2^{k-1}$. It is simply all the sequences of length $k$, composed of $0$ and $1$, and starting with a $1$. We could create this space by discovering with the algorithm all the descendants from the initial state $s_0$. 
	But in order to parallelize the algorithm, it is more convenient, and more natural, to list all the possible states at the beginning of the algorithm.
	
	We proceeded by dichotomy to find the best lower bound, that is the largest value of $p$ such that the largest eigenvalue of the transition matrix is strictly less than $1$. In order not to redo the entire algorithm for each $p$, we kept in memory the polynomials of the transitions. There are far less different polynomials than the number of states, so we listed all the different polynomials encountered, and created a first transition matrix with only the indexes of the polynomials. Then for a given value of $p$, we created the "real" transition matrix, with the corresponding values of the polynomials, a very fast procedure. 
	
	We have to apply the algorithm to a $k$ as large as possible to get the best possible lower bound. But when $k$ is changed to $k+1$, the number of states is doubled, and the number of transitions is multiplied by four, a far too abrupt change. 
	As we wanted to use all the available memory of our computer, we took a space characterized by three values $k$, $i$ and $j$, and denoted $S(k, i, j)$. This space is a composed of sequences of $0$ and $1$ of length $k+1$ starting with $1$, with
	\begin{itemize}
		\item all the sequences ending with $0$ 
		\item the sequences ending with $1$ and having at most $i+1$ ones
		\item exactly $j$ sequences ending with $1$ and having exactly $i+2$ ones.
	\end{itemize}
	Hence $S(k, i, j)$ contains all the previous sequences of length $k$ and is contained in the sequences of length $k+1$, and so can be used as an intermediate between these two spaces.
	When a potential child of a state does not belong to $S(k, i, j)$, we simply replace the one present at the last position (we know that it is the case) by a zero.
	The results indicated in section $1$ for the bond and site oriented percolation were obtained with the space $S(16,7,3620)$, with $46337$ states. Moreover, this allowed us to get precise comparisons with the results of \cite{belitsky}, see below.

	\subsection{Comparison with the method of \cite{belitsky}} \label{sec:comparaison}
	In \cite{belitsky}, the authors obtained that $p_c\geq 0.6338$ for the oriented bond percolation model. This is the most recent article we are aware of concerning a lower bound in the classical two-dimensional oriented percolation, and their method could have been applied in the site percolation setting.

	In an attempt to give a fair comparison between this result and ours (that is $p_c\geq 0.636893$), we list their results for different sizes they considered, together with the number of states we needed to obtain equivalent bounds.
	
	\begin{center}
		\begin{tabular}{ |l|r|p{20mm}|p{28mm}|r| } 
			\hline
			Lower bound in \cite{belitsky}& $k$ in  \cite{belitsky} &number of states in \cite{belitsky} & present number of states, and space & ratio \\ 
			\hline
			$0.624211$ &$4$ &$16$ & $38$ $S(6,2,0)$ &$2.375$ \\ 
			$0.627067$ & $5$&$32$ & $71$ $S(7,2,0)$& $2.2188$\\
			$0.629203$ & $6$ & $64$ & $137$ $S(8,2,1)$ & $2.1407$ \\
			$0.630864$ & $7$ &$128$  & $275$ $S(9,2,10)$ & $2.1485$ \\
			$0.632193$  & $8$ &$256$  & $550$ $S(10,2,28)$ & $2.1485$ \\
			$0.63328$   & $9$ & $512$ & $1073$ $S(11,2,44)$&  $2.0958$\\
			\hline
		\end{tabular}
	\end{center}
	The ratio is the ratio of our number of states to theirs. For example, in the first line, we had to use a space of size $38$ to get a lower bound at least as good as $0.624211$, the ratio is obtained by $2.375=38/16$.
	
	It may be conjectured that we obtain similar results with a state space approximately $2.1$ times that of \cite{belitsky}, a comparaison which is not in favor of our new mehod.
	But there were in fact two transition matrices in \cite{belitsky}:
	\begin{itemize}
		
		\item $plm$ for the transition matrix between different segments of size $k$ starting from the rightmost infected vertex.
		\item $plmk$  for the transition matrix between different segments of size $k$, starting from the rightmost infected vertex, \underline{and} the offset between two consecutives layers.
	\end{itemize}
	
	On one hand, the matrix of $plm$ is indeed approximately $4.2$ times smaller than our transition matrix. They used their matrix to calculate an invariant measure, while we approximate the largest eigenvalue of our matrix.
	
	On the other hand, the matrix $plmk$ is about $k$ times larger than $plm$, and for $k=14$, this gives more than $3$ times the size of our matrix, a value which grows linearly with $k$.
	
	So with their method, one works (finding an invariant measure or the largest eigenvalue) with a smaller matrix (more than four times), but one has to keep in memory three times more coefficients.
	
	Another discrepancy between the two methods lies in the fact that they used a result concerning the rightmost infected particle, based on a result due to Durrett~\cite{Durrett}. 
	In the next section, we shall present variants of the two-dimensional oriented percolation. While the result concerning the connection between the rightmost particle and the existence of percolation is certainly provable for the second variant, our method can be applied with a simple change of the transition probabilities. 
	
	One way to circumvent this issue concerning the method of \cite{belitsky} could be of course to study the relative position of the rightmost particle and the leftmost one, which should work on all cases, but with a twofold impact on the memory.
	
	Finally, we are going to highlight in section~\ref{sec:threeDim} that our method does not rely on topological considerations, and apply it in the three-dimensional setting. We will not pursue the cases of larger dimensions, but it would be nevertheless applicable.
	
	\subsection{Alternative two-dimensional graph}
	Another classical two-dimensional oriented graph is $\Laltd$, where from each vertex $(x, y)$ we have three oriented edges toward $(x-1, y+1)$, $(x, y+1)$ and $(x+1, y+1)$, see below. 
	
	
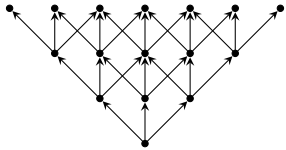
\begin{center}
		\begin{tikzpicture}[scale=0.6, ->]
			\filldraw (0,0) circle (2pt);
			\filldraw (1,1) circle (2pt);
			\filldraw (-1,1) circle (2pt);
			\filldraw (0,1) circle (2pt);
			\filldraw (1,2) circle (2pt);
			\filldraw (-1,2) circle (2pt);
			\filldraw (0,2) circle (2pt);
			\filldraw (2,2) circle (2pt);
			\filldraw (-2,2) circle (2pt);
			\filldraw (1,3) circle (2pt);
			\filldraw (-1,3) circle (2pt);
			\filldraw (0,3) circle (2pt);
			\filldraw (2,3) circle (2pt);
			\filldraw (-2,3) circle (2pt);
			\filldraw (3,3) circle (2pt);
			\filldraw (-3,3) circle (2pt);
			\draw[>=stealth, shorten >= 1.5pt](0,0)--(0,1);
			\draw[>=stealth, shorten >= 1.5pt](0,0)--(-1,1);
			\draw[>=stealth, shorten >= 1.5pt](0,0)--(1,1);
			\begin{scope}[shift={(0,1)}]
				\draw[>=stealth, shorten >= 1.5pt](0,0)--(0,1);
				\draw[>=stealth, shorten >= 1.5pt](0,0)--(-1,1);
				\draw[>=stealth, shorten >= 1.5pt](0,0)--(1,1);
			\end{scope}
			\begin{scope}[shift={(-1,1)}]
				\draw[>=stealth, shorten >= 1.5pt](0,0)--(0,1);
				\draw[>=stealth, shorten >= 1.5pt](0,0)--(-1,1);
				\draw[>=stealth, shorten >= 1.5pt](0,0)--(1,1);
			\end{scope}
			\begin{scope}[shift={(1,1)}]
				\draw[>=stealth, shorten >= 1.5pt](0,0)--(0,1);
				\draw[>=stealth, shorten >= 1.5pt](0,0)--(-1,1);
				\draw[>=stealth, shorten >= 1.5pt](0,0)--(1,1);
			\end{scope}
			\begin{scope}[shift={(1,2)}]
				\draw[>=stealth, shorten >= 1.5pt](0,0)--(0,1);
				\draw[>=stealth, shorten >= 1.5pt](0,0)--(-1,1);
				\draw[>=stealth, shorten >= 1.5pt](0,0)--(1,1);
			\end{scope}
			\begin{scope}[shift={(2,2)}]
				\draw[>=stealth, shorten >= 1.5pt](0,0)--(0,1);
				\draw[>=stealth, shorten >= 1.5pt](0,0)--(-1,1);
				\draw[>=stealth, shorten >= 1.5pt](0,0)--(1,1);
			\end{scope}
			\begin{scope}[shift={(0,2)}]
				\draw[>=stealth, shorten >= 1.5pt](0,0)--(0,1);
				\draw[>=stealth, shorten >= 1.5pt](0,0)--(-1,1);
				\draw[>=stealth, shorten >= 1.5pt](0,0)--(1,1);
			\end{scope}
			\begin{scope}[shift={(-1,2)}]
				\draw[>=stealth, shorten >= 1.5pt](0,0)--(0,1);
				\draw[>=stealth, shorten >= 1.5pt](0,0)--(-1,1);
				\draw[>=stealth, shorten >= 1.5pt](0,0)--(1,1);
			\end{scope}
			\begin{scope}[shift={(-2,2)}]
				\draw[>=stealth, shorten >= 1.5pt](0,0)--(0,1);
				\draw[>=stealth, shorten >= 1.5pt](0,0)--(-1,1);
				\draw[>=stealth, shorten >= 1.5pt](0,0)--(1,1);
			\end{scope}
		\end{tikzpicture}
	\end{center}
	
	An edge is said {\it good} if the infected vertex is not also infected from a vertex on the right.
	
	\begin{center}
		\begin{tikzpicture}[scale=0.6, ->]
			\filldraw (0,0) circle (2pt);
			\filldraw (1,1) circle (2pt);
			\filldraw (-1,1) circle (2pt);
			\filldraw (0,1) circle (2pt);
			\filldraw (1,2) circle (2pt);
			\draw (-1,2) circle (2pt);
			\filldraw (0,2) circle (2pt);
			\filldraw (2,2) circle (2pt);
			\filldraw (-2,2) circle (2pt);
			\filldraw (1,3) circle (2pt);
			\filldraw (-1,3) circle (2pt);
			\filldraw (0,3) circle (2pt);
			\filldraw (2,3) circle (2pt);
			\filldraw (-2,3) circle (2pt);
			\filldraw (3,3) circle (2pt);
			\draw (-3,3) circle (2pt);
			\draw[>=stealth, shorten >= 1.5pt](0,0)--(0,1);
			\draw[>=stealth, shorten >= 1.5pt](0,0)--(-1,1);
			\draw[>=stealth, shorten >= 1.5pt](0,0)--(1,1);
			\begin{scope}[shift={(0,1)}]
				\draw[>=stealth, shorten >= 1.5pt](0,0)--(0,1);
			\end{scope}
			\begin{scope}[shift={(-1,1)}]
				\draw[>=stealth, shorten >= 1.5pt](0,0)--(-1,1);
				\draw[>=stealth, shorten >= 1.5pt](0,0)--(1,1);
			\end{scope}
			\begin{scope}[shift={(1,1)}]
				\draw[>=stealth, shorten >= 1.5pt](0,0)--(0,1);
				\draw[>=stealth, shorten >= 1.5pt](0,0)--(-1,1);
				\draw[>=stealth, shorten >= 1.5pt](0,0)--(1,1);
			\end{scope}
			\begin{scope}[shift={(1,2)}]
				\draw[>=stealth, shorten >= 1.5pt](0,0)--(-1,1);
				\draw[>=stealth, shorten >= 1.5pt](0,0)--(1,1);
			\end{scope}
			\begin{scope}[shift={(2,2)}]
				\draw[>=stealth, shorten >= 1.5pt](0,0)--(-1,1);
				\draw[>=stealth, shorten >= 1.5pt](0,0)--(1,1);
			\end{scope}
			\begin{scope}[shift={(0,2)}]
				\draw[>=stealth, shorten >= 1.5pt](0,0)--(0,1);
				\draw[>=stealth, shorten >= 1.5pt](0,0)--(-1,1);
				\draw[>=stealth, shorten >= 1.5pt](0,0)--(1,1);
			\end{scope}
			\begin{scope}[shift={(-1,2)}]
				\draw[>=stealth, shorten >= 1.5pt](0,0)--(-1,1);
				\draw[>=stealth, shorten >= 1.5pt](0,0)--(1,1);
			\end{scope}
			\begin{scope}[shift={(-2,2)}]
				\draw[>=stealth, shorten >= 1.5pt](0,0)--(0,1);
				\draw[>=stealth, shorten >= 1.5pt](0,0)--(1,1);
			\end{scope}
			\node at (0,-0.8) {open edges};
			\begin{scope}[shift={(8.5,0)}]
				\node at (0,-0.8) {good edges};
				\filldraw (0,0) circle (2pt);
				\filldraw (1,1) circle (2pt);
				\filldraw (-1,1) circle (2pt);
				\filldraw (0,1) circle (2pt);
				\filldraw (1,2) circle (2pt);
				\draw (-1,2) circle (2pt);
				\filldraw (0,2) circle (2pt);
				\filldraw (2,2) circle (2pt);
				\filldraw (-2,2) circle (2pt);
				\filldraw (1,3) circle (2pt);
				\filldraw (-1,3) circle (2pt);
				\filldraw (0,3) circle (2pt);
				\filldraw (2,3) circle (2pt);
				\filldraw (-2,3) circle (2pt);
				\filldraw (3,3) circle (2pt);
				\draw (-3,3) circle (2pt);
				\draw[>=stealth, shorten >= 1.5pt](0,0)--(0,1);
				\draw[>=stealth, shorten >= 1.5pt](0,0)--(-1,1);
				\draw[>=stealth, shorten >= 1.5pt](0,0)--(1,1);
				\begin{scope}[shift={(0,1)}]
				\end{scope}
				\begin{scope}[shift={(-1,1)}]
					\draw[>=stealth, shorten >= 1.5pt](0,0)--(-1,1);
				\end{scope}
				\begin{scope}[shift={(1,1)}]
					\draw[>=stealth, shorten >= 1.5pt](0,0)--(0,1);
					\draw[>=stealth, shorten >= 1.5pt](0,0)--(-1,1);
					\draw[>=stealth, shorten >= 1.5pt](0,0)--(1,1);
				\end{scope}
				\begin{scope}[shift={(1,2)}]
					\draw[>=stealth, shorten >= 1.5pt](0,0)--(-1,1);
					\draw[>=stealth, shorten >= 1.5pt](0,0)--(1,1);
				\end{scope}
				\begin{scope}[shift={(2,2)}]
					\draw[>=stealth, shorten >= 1.5pt](0,0)--(-1,1);
					\draw[>=stealth, shorten >= 1.5pt](0,0)--(1,1);
				\end{scope}
				\begin{scope}[shift={(0,2)}]
					\draw[>=stealth, shorten >= 1.5pt](0,0)--(-1,1);
				\end{scope}
				\begin{scope}[shift={(-2,2)}]
					\draw[>=stealth, shorten >= 1.5pt](0,0)--(0,1);
				\end{scope}
			\end{scope}
		\end{tikzpicture}
	\end{center}

	By a Peierl argument, we know that the critical points, both for the site and the bond percolation, are larger than $1/3$.
	In \cite{Pearce}, it was proved furthermore that, for site percolation, $p_c\geq 0.5$.
	Here we obtained the lower bounds $0.525$ with $k=15$ for site percolation, and $0.4022$ with $k=13$ for bond percolation.
	
	We used only $k=13$ for bond percolation since, in this case, the transition probabilities are polynomials of somewhat high degrees.
	
	
	\section{Inhomogeneous two-dimensional oriented percolation processes}
	\label{sec:otherTwoDim}
	We show here that our technique is adaptable to various inhomogeneous cases on $\Z^2$.
	\subsection{Presentation of the models}
	All the models here have two parameters $p_1$ and $p_2$ in $[0, 1]$. The parameter $p_2$ is considered fixed, and we seek a lower bound for the critical value of $p_1$. The variables $x, y, n$ are integers.
	\begin{itemize}
		\item Model I: this is a vertex percolation model. Each vertex of the form $(2n, y)$ is open with probability $p_1$, while each vertex of the form $(2n+1, y)$ is open with probability $p_2$.
		\item Model II: a vertex percolation model. Each vertex of the form $(2n-y, y)$ is open with probability $p_1$, and each vertex of the form $(2n+1-y, y)$ is open with probability $p_2$.
		\item Model III: a vertex percolation model. Consider a vertex $(x, y)$. If $\mod(x-y, 4)\leq 1$, then the vertex is open with probability $p_1$. Otherwise it is open with probability $p_2$.
		\item Model IV: a vertex percolation model. Consider a vertex $(x, y)$. If $\mod(x, 4)\leq 1$, then the vertex is open with probability $p_1$. Otherwise it is open with probability $p_2$.
		\item Model V: an edge percolation model. Horizontal edges are open with probability $p_1$, and vertical ones are open with probability $p_2$.
	\end{itemize}
	
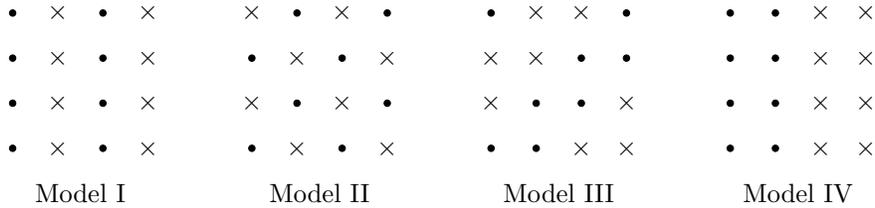
\begin{center}
		\begin{tikzpicture}[scale=0.6]
			\filldraw (0,3) circle (2pt);
			\filldraw (0,2) circle (2pt);
			\filldraw (0,1) circle (2pt);
			\filldraw (0,0) circle (2pt);
			\node at(1,0) {$\times$};
			\node at(1,1) {$\times$};
			\node at(1,2) {$\times$};
			\node at(1,3) {$\times$};
			\begin{scope}[shift={(2,0)}]
				\filldraw (0,3) circle (2pt);
				\filldraw (0,2) circle (2pt);
				\filldraw (0,1) circle (2pt);
				\filldraw (0,0) circle (2pt);
				\node at(1,0) {$\times$};
				\node at(1,1) {$\times$};
				\node at(1,2) {$\times$};
				\node at(1,3) {$\times$};
			\end{scope}
			\node at (1.5, -1) {Model I};
			\begin{scope}[shift={(5.3,0)}]
				\filldraw (3,1) circle (2pt);
				\filldraw (2,0) circle (2pt);
				\filldraw (1,1) circle (2pt);
				\filldraw (0,0) circle (2pt);
				\filldraw (0,2) circle (2pt);
				\filldraw (2,2) circle (2pt);
				\filldraw (1,3) circle (2pt);
				\filldraw (3,3) circle (2pt);
				\node at(1,0) {$\times$};
				\node at(0,1) {$\times$};
				\node at(1,2) {$\times$};
				\node at(2,1) {$\times$};
				\node at(3,0) {$\times$};
				\node at(0,3) {$\times$};
				\node at(3,2) {$\times$};
				\node at(2,3) {$\times$};
				\node at (1.5, -1) {Model II};
			\end{scope}
			\begin{scope}[shift={(10.6,0)}]
				\filldraw (2,1) circle (2pt);
				\filldraw (1,1) circle (2pt);
				\filldraw (1,0) circle (2pt);
				\filldraw (0,0) circle (2pt);
				\filldraw (3,2) circle (2pt);
				\filldraw (2,2) circle (2pt);
				\filldraw (3,3) circle (2pt);
				\filldraw (0,3) circle (2pt);
				\node at(0,2) {$\times$};
				\node at(0,1) {$\times$};
				\node at(1,2) {$\times$};
				\node at(1,3) {$\times$};
				\node at(2,0) {$\times$};
				\node at(3,0) {$\times$};
				\node at(3,1) {$\times$};
				\node at(2,3) {$\times$};
				\node at (1.5, -1) {Model III};
			\end{scope}
			\begin{scope}[shift={(15.9,0)}]
				\filldraw (0,3) circle (2pt);
				\filldraw (0,2) circle (2pt);
				\filldraw (0,1) circle (2pt);
				\filldraw (0,0) circle (2pt);
				\node at(2,0) {$\times$};
				\node at(2,1) {$\times$};
				\node at(2,2) {$\times$};
				\node at(2,3) {$\times$};
				\begin{scope}[shift={(1,0)}]
					\filldraw (0,3) circle (2pt);
					\filldraw (0,2) circle (2pt);
					\filldraw (0,1) circle (2pt);
					\filldraw (0,0) circle (2pt);
					\node at(2,0) {$\times$};
					\node at(2,1) {$\times$};
					\node at(2,2) {$\times$};
					\node at(2,3) {$\times$};
				\end{scope}
				\node at (1.5, -1) {Model IV};
			\end{scope}
			
			\filldraw (2,-2.6) circle (2pt);
			\node at (8, -2.6) { are associated to $p_1$, $\times$ are associated to $p_2$};
		\end{tikzpicture}
	\end{center}
	
	\subsection{The Lower Bounds}
	We report in the next table the lower bounds for the critical value of $p_1$ for the different models and various values of $p_2$. Note that the sizes of the automatons, to be described after, depend on the models, and so we used different values for $k$.
	
	\begin{center}
		\begin{tabular}{ |l|l|c|c| } 
			\hline
			Model & $p_2$ &Lower bound & Value of $k$ used \\ \hline
			Model I & $0.6$&$0.7693$ & $15$ \\ 
			Model I & $0.8$&$0.5444$ & $15$ \\ 
			Model II & $0.6$&$0.8189$ & $15$ \\ 
			Model II & $0.8$&$0.6103$ & $14$ \\
			Model II & $1$&$0.5223$ & $14$ \\
			Model III & $0.6$&$0.7759$ & $13$ \\ 
			Model III & $0.8$& $0.5753$& $13$ \\ 
			Model IV & $0.6$&$0.7720$ & $13$ \\ 
			Model IV & $0.8$&$0.5583$ & $13$ \\ 
			Model V & $0.5$&$0.7539$ & $15$ \\ 
			\hline
		\end{tabular}
	\end{center}
	Note that for Model II, the critical value of $p_1$ is not trivial even when $p_2=1$.
	For the vertex models, the lower bounds are in the following order:
	
	\[\mbox{Model I} < \mbox{Model IV} < \mbox{Model III} < \mbox{Model II}.\]
	We conjecture that the critical values verify the same order.
	\subsection{The automatons}
	We shall present the details for Model I, and draw upon these explanations for the other models.
	\subsubsection{Model I}
	As before a state is associated to a segment of length $k$. A vertex is said of type $A$ if it is associated to $p_1$, and of type $B$ if it is associated to $p_2$. Now a state is the sequence of the status of the vertices of the segment (occupied or vacant), {\it together } with a function of the coordinates $(x, y)$ of its rightmost vertex. Here, if the rightmost vertex is of type $A$, then this function returns $a$, and in the other case it returns $b$. For example, 
	
	\begin{center}
		\begin{tikzpicture}[scale=0.6]
			\filldraw (0,0) circle (2pt);
			\filldraw (2,0) circle (2pt);
			\node at(1,0) {$\times$};
			\node at(3,0) {$\times$};
			\node at(0,0.5) {$1$};
			\node at(1,0.5) {$1$};
			\node at(2,0.5) {$0$};
			\node at(3,0.5) {$0$};
			\node at(2.5,-0.9) {is the state $(1100,b)$};
			\begin{scope}[shift={(8.2,0)}]
				\filldraw (1,0) circle (2pt);
				\filldraw (3,0) circle (2pt);
				\node at(0,0) {$\times$};
				\node at(2,0) {$\times$};
				\node at(0,0.5) {$1$};
				\node at(1,0.5) {$0$};
				\node at(2,0.5) {$0$};
				\node at(3,0.5) {$1$};
				\node at(2.5,-0.8) {is the state $(1001,a)$};
			\end{scope}
		\end{tikzpicture}
	\end{center}
	In Model I, the types of the vertices are alternating. On the following picture:
	
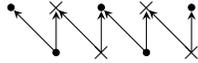
\begin{center}
		\begin{tikzpicture}[scale=0.6, ->]
			\draw[>=stealth, shorten >= 1.5pt] (0,0)--(-1,1);
			\draw[>=stealth, shorten >= 1.5pt] (0,0)--(0,1);
			\draw[>=stealth, shorten >= 1.5pt] (1,0)--(0,1);
			\draw[>=stealth, shorten >= 1.5pt] (1,0)--(1,1);
			\draw[>=stealth, shorten >= 1.5pt] (2,0)--(1,1);
			\draw[>=stealth, shorten >= 1.5pt] (2,0)--(2,1);
			\draw[>=stealth, shorten >= 1.5pt] (3,0)--(2,1);
			\draw[>=stealth, shorten >= 1.5pt] (3,0)--(3,1);
			\filldraw (0,0) circle (2pt);
			\filldraw (2,0) circle (2pt);
			\node at(1,0) {$\times$};
			\node at(3,0) {$\times$};
			\begin{scope}[shift={(0,1)}]
				\filldraw (-1,0) circle (2pt);
				\filldraw (1,0) circle (2pt);
				\filldraw (3,0) circle (2pt);
				\node at(0,0) {$\times$};
				\node at(2,0) {$\times$};
			\end{scope}
		\end{tikzpicture}
	\end{center}
	one can see that the up-child has the same variant as its parent, and the right-child has the other variant.
	\subsubsection{Model II}
	In Model II, the types of the vertices in a state are all identical. We keep the same denomination for the states, and adapt the transitions.
	
	This time, both the up-child and the right-child have the other variant than their parent.
	\subsubsection{Model III}
	In Model III, the types of the vertices are alternating, as in Model I, but the function associated to the denomination of a state will now return four different values:
	\begin{itemize}
		\item if the rightmost vertex is of type $A$ and its successor on the right is also of type $A$, then the function returns $a$
		\item if the rightmost vertex is of type $A$ and its successor on the right is  of type $B$, then the function returns $b$
		\item if the rightmost vertex is of type $B$ and its successor on the right is of type $A$, then the function returns $c$
		\item if the rightmost vertex is of type $B$ and its successor on the right is also of type $B$, then the function returns $d$
	\end{itemize}
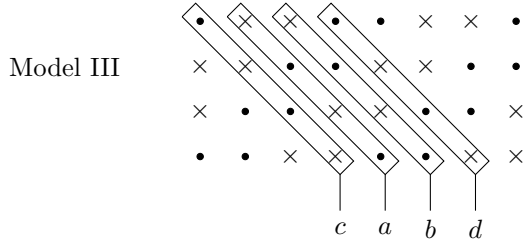
	\begin{center}
		\begin{tikzpicture}[scale=0.6]
			\draw (-0.4,3.1)--(3.1,-0.4)--(3.4,-0.1)--(-0.1,3.4)--(-0.4,3.1);
			\draw (3.1, -0.4)--(3.1, -1.2);
			\node at(3.1, -1.6) {$c$};
			\begin{scope}[shift={(1,0)}]
				\draw (-0.4,3.1)--(3.1,-0.4)--(3.4,-0.1)--(-0.1,3.4)--(-0.4,3.1);
				\draw (3.1, -0.4)--(3.1, -1.2);
				\node at(3.1, -1.6) {$a$};
			\end{scope}
			\begin{scope}[shift={(2,0)}]
				\draw (-0.4,3.1)--(3.1,-0.4)--(3.4,-0.1)--(-0.1,3.4)--(-0.4,3.1);
				\draw (3.1, -0.4)--(3.1, -1.2);
				\node at(3.1, -1.6) {$b$};
			\end{scope}
			\begin{scope}[shift={(3,0)}]
				\draw (-0.4,3.1)--(3.1,-0.4)--(3.4,-0.1)--(-0.1,3.4)--(-0.4,3.1);
				\draw (3.1, -0.4)--(3.1, -1.2);
				\node at(3.1, -1.6) {$d$};
			\end{scope}
			\filldraw (2,1) circle (2pt);
			\filldraw (1,1) circle (2pt);
			\filldraw (1,0) circle (2pt);
			\filldraw (0,0) circle (2pt);
			\filldraw (3,2) circle (2pt);
			\filldraw (2,2) circle (2pt);
			\filldraw (3,3) circle (2pt);
			\filldraw (0,3) circle (2pt);
			\node at(0,2) {$\times$};
			\node at(0,1) {$\times$};
			\node at(1,2) {$\times$};
			\node at(1,3) {$\times$};
			\node at(2,0) {$\times$};
			\node at(3,0) {$\times$};
			\node at(3,1) {$\times$};
			\node at(2,3) {$\times$};
			\begin{scope}[shift={(4,0)}]
				\filldraw (2,1) circle (2pt);
				\filldraw (1,1) circle (2pt);
				\filldraw (1,0) circle (2pt);
				\filldraw (0,0) circle (2pt);
				\filldraw (3,2) circle (2pt);
				\filldraw (2,2) circle (2pt);
				\filldraw (3,3) circle (2pt);
				\filldraw (0,3) circle (2pt);
				\node at(0,2) {$\times$};
				\node at(0,1) {$\times$};
				\node at(1,2) {$\times$};
				\node at(1,3) {$\times$};
				\node at(2,0) {$\times$};
				\node at(3,0) {$\times$};
				\node at(3,1) {$\times$};
				\node at(2,3) {$\times$};
			\end{scope}
			\node at (-3, 2) {Model III};
		\end{tikzpicture}
	\end{center}
	The variants of the children are as follows :
	\begin{center}
		\begin{tabular}{ |c|c|c| } 
			\hline
			Parent & up-child & right-child \\ 
			\hline
			$a$ & $c$ &$b$ \\ 
			$b$ & $a$ & $d$ \\ 
			$c$& $d$& $a$\\
			$d$ & $b$&$c$\\ 
			\hline
		\end{tabular}
	\end{center}
	
	\subsubsection{Model IV}
	For Model IV, we keep the same function as in Model III.
	
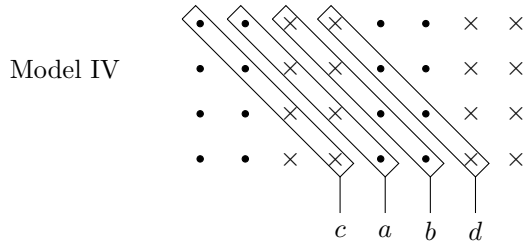
\begin{center}
		\begin{tikzpicture}[scale=0.6]
			\draw (-0.4,3.1)--(3.1,-0.4)--(3.4,-0.1)--(-0.1,3.4)--(-0.4,3.1);
			\draw (3.1, -0.4)--(3.1, -1.2);
			\node at(3.1, -1.6) {$c$};
			\begin{scope}[shift={(1,0)}]
				\draw (-0.4,3.1)--(3.1,-0.4)--(3.4,-0.1)--(-0.1,3.4)--(-0.4,3.1);
				\draw (3.1, -0.4)--(3.1, -1.2);
				\node at(3.1, -1.6) {$a$};
			\end{scope}
			\begin{scope}[shift={(2,0)}]
				\draw (-0.4,3.1)--(3.1,-0.4)--(3.4,-0.1)--(-0.1,3.4)--(-0.4,3.1);
				\draw (3.1, -0.4)--(3.1, -1.2);
				\node at(3.1, -1.6) {$b$};
			\end{scope}
			\begin{scope}[shift={(3,0)}]
				\draw (-0.4,3.1)--(3.1,-0.4)--(3.4,-0.1)--(-0.1,3.4)--(-0.4,3.1);
				\draw (3.1, -0.4)--(3.1, -1.2);
				\node at(3.1, -1.6) {$d$};
			\end{scope}
			\filldraw (0,3) circle (2pt);
			\filldraw (0,2) circle (2pt);
			\filldraw (0,1) circle (2pt);
			\filldraw (0,0) circle (2pt);
			\node at(2,0) {$\times$};
			\node at(2,1) {$\times$};
			\node at(2,2) {$\times$};
			\node at(2,3) {$\times$};
			\begin{scope}[shift={(1,0)}]
				\filldraw (0,3) circle (2pt);
				\filldraw (0,2) circle (2pt);
				\filldraw (0,1) circle (2pt);
				\filldraw (0,0) circle (2pt);
				\node at(2,0) {$\times$};
				\node at(2,1) {$\times$};
				\node at(2,2) {$\times$};
				\node at(2,3) {$\times$};
			\end{scope}
			\begin{scope}[shift={(4,0)}]
				\filldraw (0,3) circle (2pt);
				\filldraw (0,2) circle (2pt);
				\filldraw (0,1) circle (2pt);
				\filldraw (0,0) circle (2pt);
				\node at(2,0) {$\times$};
				\node at(2,1) {$\times$};
				\node at(2,2) {$\times$};
				\node at(2,3) {$\times$};
			\end{scope}
			\begin{scope}[shift={(5,0)}]
				\filldraw (0,3) circle (2pt);
				\filldraw (0,2) circle (2pt);
				\filldraw (0,1) circle (2pt);
				\filldraw (0,0) circle (2pt);
				\node at(2,0) {$\times$};
				\node at(2,1) {$\times$};
				\node at(2,2) {$\times$};
				\node at(2,3) {$\times$};
			\end{scope}
			\node at (-3, 2) {Model IV};
		\end{tikzpicture}
	\end{center}
	The variants of the children are in this case :
	\begin{center}
		\begin{tabular}{ |c|c|c| } 
			\hline
			Parent & up-child & right-child \\ 
			\hline
			$a$ & $a$ &$b$ \\ 
			$b$ & $b$ & $d$ \\ 
			$c$& $c$& $a$\\
			$d$ & $d$&$c$\\ 
			\hline
		\end{tabular}
	\end{center}
	\subsubsection{Model V}
	Model V does not need an auxiliary function, as all vertices are equivalent.
	\section{Three-dimensional oriented percolation}
	\label{sec:threeDim}
	We study here the classical three-dimensional oriented lattice $\Lva$, with  the site and bond percolation variants.
	
	We denote by $e_1, e_2, e_3$ the three canonical vectors in dimension $3$. In dimension $2$, the vector $e_2$ had priority over $e_1$. Here $e_3$ has priority over $e_2$, which in turns has priority over $e_1$. 
	
	Again, for the explanations, we consider only the bond percolation model. An edge, starting from the vertex $(x, y, z)$, is said {\it good} if it is open, starts from an occupied vertex, and either
	\begin{itemize}
		\item is in the direction $e_3$.
		\item is in the direction $e_2$, and there is no open path from the origin to $(x, y+1, z)$ and passing through $(x, y+1, z-1)$
		\item is in the direction $e_1$, and there is no open path from the origin to $(x+1, y, z)$ and passing through $(x+1, y, z-1)$ or $(x+1, y-1, z)$.
	\end{itemize}
	The height of a vertex $(x, y, z)$ is $x+y+z$. We represent in dimension two several vertices of the same height with the following conventions:
	\begin{itemize}
		\item When we go upward, we add $(-1, 0, 1)$.
		\item When we go to the right, we add $(-1, 1, 0)$.
	\end{itemize}
	Here is a representation of ten numbered vertices of the same height:
	\begin{center}
		\begin{tikzpicture}[scale=1.2]
			\draw (0,4)--(0,0);
			\draw (1, 4)--(1,0);
			\draw (2,3)--(2,0);
			\draw (3,2)--(3,0);
			\draw (4,1)--(4,0);
			\draw (0,4)--(1,4);
			\draw (0,3)--(2,3);
			\draw (0,2)--(3,2);
			\draw (0,1)--(4,1);
			\draw (0,0)--(4,0);
			\node at(0.5, 3.72) {$1$};
			\node at(0.5, 3.23) {$(0,0,3)$};
			\node at(0.5, 2.72) {$2$};
			\node at(0.5, 2.23) {$(1,0,2)$};
			\node at(1.5, 2.72) {$3$};
			\node at(1.5, 2.23) {$(0,1,2)$};
			\node at(0.5, 1.72) {$4$};
			\node at(0.5, 1.23) {$(2,0,1)$};
			\node at(1.5, 1.72) {$5$};
			\node at(1.5, 1.23) {$(1,1,1)$};
			\node at(2.5, 1.72) {$6$};
			\node at(2.5, 1.23) {$(0,2,1)$};
			\node at(0.5, 0.72) {$7$};
			\node at(0.5, 0.23) {$(3,0,0)$};
			\node at(1.5, 0.72) {$8$};
			\node at(1.5, 0.23) {$(2,1,0)$};
			\node at(2.5, 0.72) {$9$};
			\node at(2.5, 0.23) {$(1,2,0)$};
			\node at(3.5, 0.72) {$10$};
			\node at(3.5, 0.23) {$(0,3,0)$};
			\begin{scope}[shift={(1,0)}]
			\end{scope}
		\end{tikzpicture}
	\end{center}
	
	For our first attempt, we used ten vertices in this disposition for a state. The vertex we focus on is the vertex numbered $3$. Note that a state can have up to three children. In our example, the vertex $2$ has for successors $(1,0,3)$, $(1,1,2)$ and $(2,0,2)$. The site $(1,1,2)$ is also a successor of the vertex $3$. So, with the rules of precedence we have chosen, if the vertex $2$ is occupied, and its "$e_2$"-edge is open, then the "$e_1$"-edge from the vertex $3$ cannot be good, and the state will have at most two children. As a last remark, the site $(1,1,2)$ is a common successor of the vertices $2$, $3$, and $5$, and three vertices having a common successor are always in this configuration.
	
	To get the children of a state, we place it in the middle of a square of size six by six. Then we place the geometric figure
	
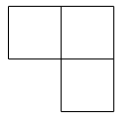
\begin{center}
		\begin{tikzpicture}[scale=0.7]
			\draw (0,2)--(0,1);
			\draw (1, 2)--(1,0);
			\draw (2,2)--(2,0);
			\draw (0,2)--(2,2);
			\draw (0,1)--(2,1);
			\draw (1,0)--(2,0);
			\begin{scope}[shift={(1,0)}]
			\end{scope}
		\end{tikzpicture}
	\end{center}
	on each position that intersects the vertices of the state. We obtain in that way the $15$ successors of the $10$ vertices, disposed in a triangle of length $5$. Then the three potential children are associated to the three triangles of length $4$ included in the triangle of length $5$.
	
	
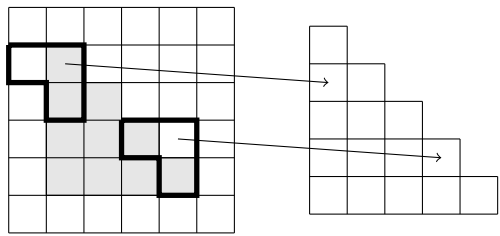
\begin{center}
		\begin{tikzpicture}[scale=0.5]
			\draw[black,fill=gray!20] (1,1)--(5,1)--(5,2)--(4,2)--(4,3)--(3,3)--(3,4)--(2,4)--(2,5)--(1,5)--(1,1);
			\draw (0,6)--(0,0);
			\draw (1, 6)--(1,0);
			\draw (2,6)--(2,0);
			\draw (3,6)--(3,0);
			\draw (4, 6)--(4,0);
			\draw (5,6)--(5,0);
			\draw (6,6)--(6,0);
			\draw (0,6)--(6,6);
			\draw (0,5)--(6,5);
			\draw (0,4)--(6,4);
			\draw (0,3)--(6,3);
			\draw (0,2)--(6,2);
			\draw (0,1)--(6,1);
			\draw (0,0)--(6,0);
			\begin{scope}[shift={(-2,1)}]
				\draw [line width=2pt] (2,4)--(4,4)--(4,2)--(3,2)--(3,3)--(2,3)--(2,4);
			\end{scope}
			\begin{scope}[shift={(1,-1)}]
				\draw [line width=2pt] (2,4)--(4,4)--(4,2)--(3,2)--(3,3)--(2,3)--(2,4);
			\end{scope}
			\draw [->] (1.5,4.5)--(8.5,4);
			\draw [->] (4.5,2.5)--(11.5,2);
			\begin{scope}[shift={(8,1.5)}]
				\draw (0,4)--(0,-1);
				\draw (1, 4)--(1,-1);
				\draw (2,3)--(2,-1);
				\draw (3,2)--(3,-1);
				\draw (4,1)--(4,-1);
				\draw (5,0)--(5,-1);
				\draw (0,4)--(1,4);
				\draw (0,3)--(2,3);
				\draw (0,2)--(3,2);
				\draw (0,1)--(4,1);
				\draw (0,0)--(5,0);
				\draw  (0,-1)--(5,-1);
				
			\end{scope}
		\end{tikzpicture}
	\end{center}
	Similarly to the two-dimensional case, a triangle of length $4$ corresponds to a child if its vertex on position $3$ is infected by the vertex on position $3$ of the initial state via a good edge.
	
	We obtained with these states $0.41$ as a lower bound for site percolation, and $0.36684$ for bond percolation. 
	As the polynomials tend to be of high degree for bond percolation, we could not go further in that case. For site percolation, we could use triangles of length $5$, that is states with $15$ vertices. Focusing on vertex $6$, we obtained $0.41507$.
	We also tried the algorithm with the focus on vertex $3$, and obtained only $0.4112$.

	\bibliographystyle{plain}
	\bibliography{biblio}
\end{document}